\input amstex

\documentstyle{amsppt}
  \magnification=1100
  \hsize=6.2truein
  \vsize=9.0truein
  \hoffset 0.1truein
  \parindent=2em

\NoBlackBoxes


\font\eusm=eusm10                   


\font\eusms=eusm7                       

\font\eusmss=eusm5                      


\newcount\theTime
\newcount\theHour
\newcount\theMinute
\newcount\theMinuteTens
\newcount\theScratch
\theTime=\number\time
\theHour=\theTime
\divide\theHour by 60
\theScratch=\theHour
\multiply\theScratch by 60
\theMinute=\theTime
\advance\theMinute by -\theScratch
\theMinuteTens=\theMinute
\divide\theMinuteTens by 10
\theScratch=\theMinuteTens
\multiply\theScratch by 10
\advance\theMinute by -\theScratch

\def\today{{\number\day\space
 \ifcase\month\or
  January\or February\or March\or April\or May\or June\or
  July\or August\or September\or October\or November\or December\fi
 \space\number\year}}

\define\Ac{{\Cal A}}

\define\Afr{{\frak A}}

\define\ah{{\hat a}}

\define\alg{{\text{alg}}}

\define\Bc{{\Cal B}}

\define\biggnm#1{
  \bigg|\bigg|#1\bigg|\bigg|}

\define\bignm#1{
  \big|\big|#1\big|\big|}

\define\Bof{B}

\define\clspan{\overline\lspan}

\define\Cpx{\bold C}

\define\dif{\text{\it d}}

\define\eqdef{{\;\overset\text{def}\to=\;}}

\define\Eto#1{E_{(\to{#1})}}

\define\Fc{{\mathchoice
     {\text{\eusm F}}
     {\text{\eusm F}}
     {\text{\eusms F}}
     {\text{\eusmss F}}}}

\define\fdim{\text{\rm fdim}\,}

\define\Feu{\Fc}

\define\fpamalg#1{{\dsize\;\operatornamewithlimits*_{#1}\;}}

\define\fpiamalg#1{{\tsize\;({*_{#1}})_{\raise-.5ex\hbox{$\ssize\iota\in I$}}}}

\define\freeF{F}

\define\freeprod#1#2{\mathchoice
     {\operatornamewithlimits{\ast}_{#1}^{#2}}
     {\raise.5ex\hbox{$\dsize\operatornamewithlimits{\ast}
      _{#1}^{#2}$}\,}
     {\text{oops!}}{\text{oops!}}}

\define\freeprodi{\mathchoice
     {\operatornamewithlimits{\ast}
      _{\iota\in I}}
     {\raise.5ex\hbox{$\dsize\operatornamewithlimits{\ast}
      _{\sssize\iota\in I}$}\,}
     {\text{oops!}}{\text{oops!}}}

\define\freeprodvni{\mathchoice
      {\operatornamewithlimits{\overline{\ast}}
       _{\iota\in I}}
      {\raise.5ex\hbox{$\dsize\operatornamewithlimits{\overline{\ast}}
       _{\sssize\iota\in I}$}\,}
      {\text{oops!}}{\text{oops!}}}

\define\Ft{{\tilde F}}

\define\Hil{{\mathchoice
     {\text{\eusm H}}
     {\text{\eusm H}}
     {\text{\eusms H}}
     {\text{\eusmss H}}}}


\define\Hilto#1{\Hil_{(\to{#1})}}

\define\Integers{\bold Z}

\define\Ints{\Integers}

\define\Lambdao{{\Lambda\oup}}

\define\ld#1{{\hbox{..}(#1)\hbox{..}}}

\define\lrnm#1{\left\|#1\right\|}

\define\lspan{\text{\rm span}@,@,@,}

\define\MvN{{\Cal M}}

\define\nm#1{\|#1\|}

\define\Nats{\Naturals}

\define\Naturals{{\bold N}}

\define\NC{{\text{\rm NC}}}

\define\NCP{{\text{\rm NCP}}}

\define\NvN{{\Cal N}}

\define\Omegao{\Omega\oup}

\define\otdts#1{\otimes_{#1}\cdots\otimes_{#1}}

\define\oup{^{\text{\rm o}}}

\define\owedge{{
     \operatorname{\raise.5ex\hbox{\text{$
     \ssize{\,\bigcirc\llap{$\ssize\wedge\,$}\,}$}}}}}

\define\owedgeo#1{{
     \underset{\raise.5ex\hbox
     {\text{$\ssize#1$}}}\to\owedge}}

\define\Pc{{\Cal P}}

\define\phit{{\tilde\phi}}

\define\PSL{{\text{\rm PSL}}}

\define\psit{{\tilde\psi}}

\define\Pto#1{{P_{(\to{#1})}}}


\define\pup#1#2{{{\vphantom{#2}}^{#1}\!{#2}}\vphantom{#2}}

\define\PvN{{\Cal P}}

\define\QED{\newline
            \line{$\hfill$\qed}\enddemo}

\define\QvN{{\Cal Q}}

\define\Reals{{\bold R}}

\define\restrict{\lower .3ex
     \hbox{\text{$|$}}}

\define\RvN{{\Cal R}}

\define\Seu{\text{\eusm S}}

\define\smd#1#2{\underset{#2}\to{#1}}

\define\smdb#1#2{\undersetbrace{#2}\to{#1}}

\define\smdbp#1#2#3{\overset{#3}\to
     {\smd{#1}{#2}}}

\define\smdbpb#1#2#3{\oversetbrace{#3}\to
     {\smdb{#1}{#2}}}

\define\smdp#1#2#3{\overset{#3}\to
     {\smd{#1}{#2}}}

\define\smdpb#1#2#3{\oversetbrace{#3}\to
     {\smd{#1}{#2}}}

\define\smp#1#2{\overset{#2}\to
     {#1}}

\define\tocdots
  {\leaders\hbox to 1em{\hss.\hss}\hfill}

\define\tr{\text{\rm tr}}

\define\Tr{\text{\rm Tr}}

\define\Zt{{\tilde Z}}


  \newcount\mycitestyle \mycitestyle=1 

  \newcount\bibno \bibno=0
  \def\newbib#1{\advance\bibno by 1 \edef#1{\number\bibno}}
  \ifnum\mycitestyle=1 \def\cite#1{{\rm[\bf #1\rm]}} \fi
  \def\scite#1#2{{\rm[\bf #1\rm, #2]}}


  \newcount\ignorsec \ignorsec=0
  \def\notasec{\ignorsec=1}

  \newcount\secno \secno=0
  \def\newsec#1{\procno=0 \subsecno=0 \ignorsec=0
    \advance\secno by 1 \edef#1{\number\secno}
    \edef\currentsec{\number\secno}}

  \newcount\subsecno
  \def\newsubsec#1{\procno=0 \advance\subsecno by 1
    \edef\currentsec{\number\secno.\number\subsecno}
     \edef#1{\currentsec}}

  \newcount\appendixno \appendixno=0
  \def\newappendix#1{\procno=0 \ignorsec=0 \advance\appendixno by 1
    \ifnum\appendixno=1 \edef\appendixalpha{\hbox{A}}
      \else \ifnum\appendixno=2 \edef\appendixalpha{\hbox{B}} \fi
      \else \ifnum\appendixno=3 \edef\appendixalpha{\hbox{C}} \fi
      \else \ifnum\appendixno=4 \edef\appendixalpha{\hbox{D}} \fi
      \else \ifnum\appendixno=5 \edef\appendixalpha{\hbox{E}} \fi
      \else \ifnum\appendixno=6 \edef\appendixalpha{\hbox{F}} \fi
    \fi
    \edef#1{\appendixalpha}
    \edef\currentsec{\appendixalpha}}

  \newcount\procno \procno=0
  \def\newproc#1{\advance\procno by 1
   \ifnum\ignorsec=0 \edef#1{\currentsec.\number\procno}
                     \edef\currentproc{\currentsec.\number\procno}
   \else \edef#1{\number\procno}
         \edef\currentproc{\number\procno}
   \fi}

  \newcount\subprocno \subprocno=0
  \def\newsubproc#1{\advance\subprocno by 1
   \ifnum\subprocno=1 \edef#1{\currentproc a} \fi
   \ifnum\subprocno=2 \edef#1{\currentproc b} \fi
   \ifnum\subprocno=3 \edef#1{\currentproc c} \fi
   \ifnum\subprocno=4 \edef#1{\currentproc d} \fi
   \ifnum\subprocno=5 \edef#1{\currentproc e} \fi
   \ifnum\subprocno=6 \edef#1{\currentproc f} \fi
   \ifnum\subprocno=7 \edef#1{\currentproc g} \fi
   \ifnum\subprocno=8 \edef#1{\currentproc h} \fi
   \ifnum\subprocno=9 \edef#1{\currentproc i} \fi
   \ifnum\subprocno>9 \edef#1{TOO MANY SUBPROCS} \fi
  }

  \newcount\tagno \tagno=0
  \def\newtag#1{\advance\tagno by 1 \edef#1{\number\tagno}}



\notasec
  \newtag{\explfreeprod}
  \newtag{\Shlyred}
  \newtag{\ourred}
  \newtag{\Afreeprod}
  \newtag{\LambdaoExpl}
\newsec{\secCompr}
 \newproc{\Comproneonn}
  \newtag{\bigLo}
  \newtag{\endLo}
 \newproc{\Comprirrat}
  \newtag{\PvNdetails}
  \newtag{\freefam}
  \newtag{\inkertau}
  \newtag{\Bpic}
 \newproc{\Stickout}
  \newtag{\StickoutEq}
 \newproc{\OnlyFreeDim}
  \newtag{\eqfdim}
 \newproc{\Compr}
  \newtag{\cutIfinite}
  \newtag{\stable}
  \newtag{\cutIinfinite}
 \newproc{\FGpinfinite}
\newsec{\secAutom}
 \newproc{\fpSX}
  \newtag{\phitp}
 \newproc{\IIinf}
  \newtag{\psipexpr}
  \newtag{\Afp}
  \newtag{\matrixunits}
  \newtag{\generators}
  \newtag{\Feuone}
  \newtag{\Feutwo}
  \newtag{\Feuthr}
  \newtag{\Feufou}
  \newtag{\Feufiv}
  \newtag{\bigffam}
  \newtag{\Feutfam}
 \newproc{\IIIone}

\newbib{\AntoniLongoRad}
\newbib{\Barnett}
\newbib{\Ching}
\newbib{\ConnesZZThesis}
\newbib{\ConnesZZCountFundGp}
\newbib{\ConnesJonesZZT}
\newbib{\DykemaZZFreeProdR}
\newbib{\DykemaZZFreeDim}
\newbib{\DykemaZZInterp}
\newbib{\DykemaZZTM}
\newbib{\DykemaRordam}
\newbib{\GolodetsNessonov}
\newbib{\HaagerupZZMAP}
\newbib{\Kreweras}
\newbib{\MvNZZiv}
\newbib{\PopaZZrigid}
\newbib{\PopaZZUnivJonesAlg}
\newbib{\RadulescuZZFundGp}
\newbib{\RadulescuZZOneParamGp}
\newbib{\RadulescuZZAmalgSubfactors}
\newbib{\ShlyakhtZZamalg}
\newbib{\SpeicherZZNoncrossing}
\newbib{\SpeicherZZComb}
\newbib{\SpeicherZZFreeCal}
\newbib{\TakesakiZZDualityCrossProd}
\newbib{\VoiculescuZZSymmetries}
\newbib{\VoiculescuZZCircSemiCirc}
\newbib{\VoiculescuZZRandMat}
\newbib{\VDNbook}

\topmatter
  \title Compressions of free products of von Neumann algebras \\
         \rm (with corrections)
  \endtitle

  \author Kenneth J\. Dykema and Florin R\u adulescu \endauthor

  \date May 24, 2000 \enddate

  \rightheadtext{}

  \leftheadtext{}

  \address Dept.~of Mathematics and Computer Science,
           Odense University,
           DK-5230 Odense M, Denmark
  \endaddress

  \email dykema\@imada.ou.dk,
         \;{\it Internet}\/{\rm :}~http://www.imada.ou.dk/\~{\hskip0.1em}dykema/
  \endemail

  \address Dept.~of Mathematics,
           University of Iowa,
           Iowa City IA 52242--1466, U.S.A.
  \endaddress

  \email radulesc\@math.uiowa.edu
  \endemail

  \abstract
     A reduction formula for compressions of von Neumann algebra
     II$_1$--factors arising as free products is proved.
     This shows that the fundamental group is $\Reals^*_+$ for some such
     algebras.
     Additionally, by taking a sort of free product with an unbounded
     semicircular element, continuous one parameter groups of trace scaling
     automorphisms on II$_\infty$--factors are constructed; this produces type
     III$_1$ factors with core $\MvN\otimes B(\text{\eusms H})$, where $\MvN$
     can be a full II$_1$--factor without the Haagerup approximation property.
  \endabstract


\endtopmatter

\document \TagsOnRight \baselineskip=18pt

\heading
Introduction.
\endheading

Murray and von Neumann~\cite{\MvNZZiv} introduced for any
II$_1$--factor $\MvN$ a one parameter family of II$_1$--factors $\MvN_t$
($0<t<\infty$) associated with
$\MvN$ by a rescaling procedure which for $t<1$ consists of taking a corner
$e\MvN e$, where $e$ is a self--adjoint idempotent of trace (continuous
dimension) $t$.
The obstruction to getting new examples of II$_1$--factors in this way is the
fundamental group $\Feu(\MvN)$, consisting of all positive reals $t$ such that
$\MvN_t$ is isomorphic to $\MvN$.
Moreover, the group $\Feu(\MvN)$ is an invariant for $\MvN$.

It was only much later realized, in A\. Connes' breakthrough classification of
type III factors~\cite{\ConnesZZThesis}, that every number
$\lambda\in\Feu(\MvN)$ in the fundamental group satisfying $\lambda<1$ can be
used to construct a type III$_\lambda$ factor. 
This is done by taking an automorphism $\alpha_\lambda$ on
$\MvN\otimes B(\Hil)$ scaling the semifinite trace by $\lambda$, (which exists
if and only if $\lambda\in\Feu(\MvN)$);
then the crossed product $\bigl(\MvN\otimes B(\Hil)\bigr)\rtimes_\alpha\Ints$
yields a type III$_\lambda$ factor;
furthermore, every III$_\lambda$ factor arises in this way.
If the fundamental group of $\MvN$ is $\Reals_+^*$, the group of positive real
numbers, 
and if the associated trace scaling automorphisms can be made into a continuous
group homomorphism $t\mapsto\alpha_t$, from $\Reals_+^*$ into 
the automorphism group of $\MvN\otimes B(\Hil)$, then as
Takesaki showed in~\cite{\TakesakiZZDualityCrossProd}, a crossed product
construction,
$(\MvN\otimes B(\Hil))\rtimes_\alpha\Reals^*_+$,
yields a type III$_1$ factor;
furthermore, all III$_1$--factors arise in this way.

Murray and von Neumann~\cite{\MvNZZiv} proved that the hyperfinite
II$_1$--factor has fundamental group $\Reals_+^*$, but for other
II$_1$--factors the fundamental group remained a complete mystery until
Connes proved~\cite{\ConnesZZCountFundGp} that the II$_1$ factor of a group
with Kazhdan's property~T is countable.
Later, Golodets and Nessonov~\cite{\GolodetsNessonov} showed that every
countable subgroup of $\Reals^*_+$ is contained in the fundamental group of
some II$_1$--factor having countable fundamental group.
(See~\cite{\ConnesJonesZZT} and~\cite{\PopaZZrigid} for related results.)

In the early 80's, Voiculescu~\cite{\VoiculescuZZSymmetries} invented free
probability theory, (see also the book~\cite{\VDNbook}).
Part of this theory is the free product construction for von Neumann algebras.
Given von Neumann algebras $\Ac_\iota$, ($\iota\in I$) with normal states
$\phi_\iota$ whose GNS representations are faithful, this construction is
denoted
$$ (\Ac,\phi)=\freeprodi(\Ac_\iota,\phi_\iota), \tag{\explfreeprod} $$
where $\Ac$ is a von Neumann algebra with state $\phi$, $\Ac$ contains copies
of the $\Ac_\iota$ to which the restriction of $\phi$ is $\phi_\iota$.
An earlier version of this construction, with special conditions on the
$(\Ac_\iota,\phi_\iota)$, was devised by W.M\.~Ching~\cite{\Ching}.
In this paper we will be concerned with the free product of
II$_1$--factors $\Ac_\iota$ with respect to their unique tracial states
$\phi_\iota$ and we will abbreviate the notation~(\explfreeprod) by writing
$\Ac=\freeprodi\Ac_\iota$.
The free product of II$_1$--factors is a II$_1$--factor, as was proved first
for many cases in~\cite{\Ching}; (the general result was proved
in~\scite{\PopaZZUnivJonesAlg}{7.4} and in~\cite{\DykemaZZTM}).

Using his matrix model for freeness~\cite{\VoiculescuZZRandMat}, Voiculescu
showed~\cite{\VoiculescuZZCircSemiCirc} that
$$ L(F_n)_{1/k}\cong L(F_{1+k^2(n-1)}), $$
which implies that the fundamental group of $L(F_\infty)$ contains all positive
rational numbers.
In~\cite{\RadulescuZZFundGp}, building on Voiculescu's technique, it was shown
that the fundamental group of $L(F_\infty)$ is all of $\Reals_+^*$.
The interpolated free group factors $L(F_s)$, ($1<s\le\infty$) were found
in~\cite{\RadulescuZZAmalgSubfactors} and~\cite{\DykemaZZInterp};
thereby it was shown that either $L(F_s)\cong L(F_\infty)$ for all $s$ or
the $L(F_s)$ are all nonisomorphic and $L(F_s)$ has trivial fundamental
group for all $s<\infty$.

By considering a sort of free product of $L^\infty(\Reals)$ with an algebra
generated by an infinite semicircular element, it was
proved in~\cite{\RadulescuZZOneParamGp} that there exists a continuous one
parameter family of automorphisms of $L(F_\infty)\otimes B(\Hil)$ scaling the
trace, and hence there is a type III$_1$ factor with core
$L(F_\infty)\otimes B(\Hil)$.
For the relevance of this III$_1$ factor to quantum field theory
see~\cite{\AntoniLongoRad}.

The result about the fundamental group of $L(F_\infty)$ was generalized by
Shlyakhtenko~\cite{\ShlyakhtZZamalg}, who proved that for any II$_1$--factor
$\MvN$,
$\Feu\bigl(\MvN*L(F_\infty)\bigr)\supseteq\Feu(\MvN)$;
together with a result from~\cite{\DykemaZZFreeProdR} this gave another proof
that $\Feu(L(F_\infty))=\Reals_+^*$.
More generally and in the same paper Shlyakhtenko proved the reduction formula
for free products:
$$ \bigl(\MvN*L(F_s)\bigr)_t=(\MvN_t)*L(F_{\frac s{t^2}}), \tag{\Shlyred} $$
for every $0<t<1$.

In this paper we show that the results about $L(F_\infty)$ and
$L(F_\infty)\otimes B(\Hil)$ described above are consequences of
properties of the free product of infinitely many II$_1$--factors.
In~\S1 we prove a reduction formula for free products of II$_1$--factors,
generalizing~(\Shlyred):
if $I$ is a finite or countably infinite set and $0<t<1$ then
$$ \bigl(\freeprodi A(\iota))_t
=\bigl(\freeprodi A(\iota)_t\bigr)*L(F_{(|I|-1)(t^{-2}-1)}), \tag{\ourred} $$
whenever $(|I|-1)(t^{-2}-1)>1$, and
similar results hold when this quantity is $\le1$.
{}From this we show that if $A(n)$ is a II$_1$--factor for every $n\ge1$ and if
$$ A=\freeprod{n=1}\infty A(n) \tag{\Afreeprod} $$
is their free product then $A$ is stable under taking the free product with
$L(F_\infty)$ and the fundamental group of $A$ contains
$\bigcap_{n=1}^\infty\Feu(A(n))$.

In~\S2, by extending the technique of~\cite{\RadulescuZZOneParamGp} we define a
sort of free product involving a semifinite weight.
We use this to construct
a continuous one parameter family of trace scaling automorphisms on
$A\otimes B(\Hil)$ if $A$ is as in~(\Afreeprod) and if each
$A(n)\otimes B(\Hil)$ has a continuous one parameter family of trace scaling
automorphisms;
this gives rise to a type III$_1$ factor with core $A\otimes B(\Hil)$.
Note that it is now possible to find such examples where $A$ does not have
property $\Gamma$ of Murray and von Neumann, and is not isomorphic to
$L(F_\infty)$:
we may take one (or more) of the $A(j)$ to be $R\otimes L(G_j)$ where $R$ is
the hyperfinite II$_1$--factor and where $G_j$ is a discrete i.c.c\. group with
property~T of Kazhdan, e.g. $\PSL(m,\Ints)$ for $m\ge3$;
that a free product of II$_1$--factors is full (i.e\. non--$\Gamma$) was first
proved in many cases by Ching~\cite{\Ching}, and in full generality by
Popa~\cite{\PopaZZUnivJonesAlg} and by L\.~Barnett~\cite{\Barnett};
moreover, by an argument in~\cite{\ConnesJonesZZT}, $A$ fails to have
Haagerup's approximation property~\cite{\HaagerupZZMAP}, hence cannot be
$L(F_\infty)$.

\heading
\S\secCompr.  Compressions of free products.
\endheading

In this section we prove the formula~(\ourred) for compressions of free product
II$_1$--factors by projections.
Throughout we rely on the formulae~\cite{\DykemaZZFreeDim} for free products of
certain sorts of finite von Neumann algebras with respect to traces, which we
now summarize.
If $A$ and $B$ are finite von Neumann algebras with dimensions at least $2$ and $3$,
respectively,
and if each is finite dimensional or approximately finite dimensional or
an interpolated free group factor or a (possibly infinite) direct sum of
these, then the free product
of $A$ and $B$ with respect to given normal faithful tracial states $\tau_A$
and $\tau_B$, is either
$L(F_t)$ or $L(F_{t'})\oplus D$, where $D$ is finite dimensional;
moreover, an algorithm is given which allows one to determine whether $D$ is
present, if so to find it and $t'$, and if not to find $t$.
The algorithm to find $t$ or $t'$ above makes use of {\it free dimension},
which is a quantity $\fdim(A,\tau_A)$, assigned to every such pair
$(A,\tau_A)$.
For example, if the result of taking the free product is the factor $L(F_t)$,
then $t=\fdim(A,\tau_A)+\fdim(B,\tau_B)$.
The free dimension of finite dimensional and approximately finite dimensional
von Neumann algebras (with given tracial states) is well defined, but, for
example, we write $\fdim(L(F_s),\sigma)=s$, and we don't know if this is well
defined.
However, the only use of free dimension is to determine the parameters $t$ and
$t'$ in the above free product von Neumann algebra;
thus if the free group factors turn out to be nonisomorphic then free
dimension will be well defined in its full generality;
otherwise, the free group factors will all be isomorphic to each other and
there will be no need to determine the parameter in the first place.
In this paper, we will blithely go on using free dimension as if it were well
defined;
in any case we obtain true statements.

Let us recall some convenient notation, introduced in~\cite{\DykemaZZInterp}.
Given a family of subsets $(\Omega_k)_{k\in K}$ of an algebra, we let
$$ \Lambdao((\Omega_k)_{k\in K})=
\{x_1x_2\cdots x_m\mid m\in\Nats,\,x_j\in\Omega_{\iota_j},\,
\iota_1\ne\iota_2,\iota_2\ne\iota_3,\ldots,\iota_{m-1}\ne\iota_m\}. 
\tag{\LambdaoExpl} $$
An element $x=x_1x_2\cdots x_m$ of~(\LambdaoExpl) is called a {\it word} in
$(\Omega_k)_{k\in K}$, and the {\it letters} of this word are
$x_1,x_2,\ldots,x_m$.
Given a von Neumann algebra with a normal tracial state $\tau$,
for a subalgebra $X$ of $\MvN$ we let $X\oup$ denote $X\cap\ker\tau$.
We will furthermore use the symbol $\ominus$ to denote the orthocomplement in
$\MvN$ with respect to the inner product induced by $\tau$.

The proof of the following lemma can be viewed as a special case of the
proof of Lemma~\Comprirrat{} and, moreover, the lemma follows easily
from Theorem~1.2 of~\cite{\DykemaZZFreeDim}, ({\it viz\.} the proof of
Lemma~\Stickout);
we state it here separately from Lemma~\Comprirrat{} and write the full proof
in detail because this special
case is like the proof of Lemma~\Comprirrat{} but is more easily understood,
and because it would be slightly awkward to handle it together with
Lemma~\Comprirrat.

\proclaim{Lemma \Comproneonn}
Let $I$ be a finite or countably infinite set having at least two elements and
for each $\iota\in I$ let
$A_\iota$ be a von Neumann algebra with normal state $\phi_\iota$ whose GNS
representation is faithful.
Let $\tr_n$ denote the tracial state on $M_n(\Cpx)$ and
consider the free products of von Neumann algebras
$$ \align
(\NvN,\phi)&=\freeprodi(A_\iota,\phi_\iota) \\
(\MvN,\psi)&=\freeprodi(A_\iota\otimes M_n(\Cpx),\phi_\iota\otimes\tr_n)
\endalign $$
Fix any $\iota_0\in I$ and let
$$ p=1\otimes e\in A_{\iota_0}\otimes M_n(\Cpx)\subseteq\MvN, $$
where $e$ is a minimal projection in $M_n(\Cpx)$.
Thus $\psi(p)=1/n$.
Then
$$ (p\MvN p,n\psi\restrict_{p\MvN p})\cong(\NvN,\phi)*(L(F_k),\tau) $$
where $\tau$ is the tracial state on the free group factor $L(F_k)$ and where
$k=(|I|-1)(n^2-1)$.
\endproclaim
\demo{Proof}
Let $\pi_\iota:A_\iota\otimes M_n(\Cpx)\hookrightarrow\MvN$ be the embeddings
arising in the free product construction.
Let $p_\iota=\pi_\iota(1\otimes e)$, so that $p=p_{\iota_0}$.
Then $p$ and $p_\iota$ are free whenever
$\iota\in I\backslash\{\iota_0\}$;
it thus follows from~\scite{\DykemaZZFreeDim}{1.1} that there is a partial
isometry $v_\iota\in\{1,p,p_\iota\}''$ such that $v_\iota^*v_\iota=p$ and
$v_\iota v_\iota^*=p_\iota$;
let $v_{\iota_0}=p$.
Consider the von Neumann subalgebras
$$ \align
\QvN&=\left(\,\bigcup_{\iota\in I}
 \pi_\iota\bigl(1\otimes M_n(\Cpx)\bigr)\right)''\subseteq\MvN \\
 \vspace{1ex}
\PvN&=\left(\,\bigcup_{\iota\in I}
 v_\iota^*\pi_\iota\bigl(A_\iota\otimes e\bigr)v_\iota\right)''
 \subseteq p\MvN p.
\endalign $$
Since $\QvN$ is isomorphic to the free product with respect to traces of $|I|$
copies of $M_n(\Cpx)$, by~\cite{\DykemaZZFreeDim}
$\QvN\cong L(F_{|I|(1-n^{-2})})$.
Moreover, $\PvN$ and $\QvN$ together generate $\MvN$;
hence $\PvN$ and $p\QvN p$ together generate $p\MvN p$.
By the formula for compressions of interpolated free group
factors, $p\QvN p\cong L(F_{(|I|-1)(n^2-1)})$.

We will show that the family of $|I|+1$ algebras,
$$ p\QvN p,\,
\left(v_\iota^*\pi_\iota\bigl(A_\iota\otimes e\bigr)v_\iota\right)_{\iota\in I}
$$
is free with respect to $n\psi\restrict_{p\MvN p}$.
This will suffice to prove the lemma, since we thereby see
that $\PvN\cong\NvN$ and that $p\MvN p$ is isomorphic to the free product of
$\PvN$ and $p\QvN p$.
We must show that
$$ \Lambdao\Bigl((p\QvN p)\oup,
\bigl(v_\iota^*\pi_\iota(A_\iota\oup\otimes e)v_\iota\bigr)_{\iota\in I}\Bigr)
\subseteq\ker\psi. \tag{\bigLo} $$
Let $x$ be an element of the left hand side of~(\bigLo) which is a word of
length at least two.
Since each $v_\iota$ is in $\QvN$ (and hence is in the centralizer of $\psi$),
it follows that $v_\iota(p\QvN p)\oup v_\iota^*=(p_\iota\QvN p_\iota)\oup$.
Therefore, by peeling off the $v_\iota$ and $v_\iota^*$ from the letters in
$x$, we see that $x=x'$ for some $x'\in\Theta$, where $\Theta$ is the set of
all words
$$ a=a_1a_2\cdots a_k\in
\Lambdao\Bigl(\QvN,\,
\bigl(\pi_\iota(A_\iota\oup\otimes e)\bigr)_{\iota\in I}\Bigr)
$$
of length at least two such that if $2\le j\le n-1$, if $a_j$ comes from $\QvN$
and if both $a_{j-1}$ and $a_{j+1}$ come from
$\pi_{\iota_1}(A_{\iota_1}\oup\otimes e)$ for some
$\iota_1\in I$, then $a_j\in(p_{\iota_1}\QvN p_{\iota_1})\oup$.
Using Kaplansky's density theorem, we see that every element of $\QvN$ is the
s.o.--limit of a bounded sequence in
$$ \lspan\Bigl(\{1\}\cup
\Lambdao\bigl((\pi_\iota(1\otimes M_n(\Cpx))\oup)_{\iota\in I}\bigr)\Bigr). $$
Fix $\iota_1\in I$;
using the trace preserving conditional expectation
$E_{\iota_1}:\QvN\to\pi_{\iota_1}(1\otimes M_n(\Cpx))$ and using that
$p_{\iota_1}$
is a minimal projection in the latter algebra, we see that every element of
$(p_{\iota_1}\QvN p_{\iota_1})\oup$ is the s.o.--limit of a bounded sequence in
$$ \lspan\Bigl(\Lambdao\bigl((\pi_\iota(1\otimes M_n(\Cpx))\oup)_{\iota\in I}
\bigr)\backslash\pi_{\iota_1}(1\otimes M_n(\Cpx))\oup\Bigr), $$
namely of spans of reduced words in the $\pi_\iota(1\otimes M_n(\Cpx))\oup$ but
excluding those words consisting of single letters belonging to
$\pi_{\iota_1}(1\otimes M_n(\Cpx))\oup$.
Hence in order to show $\Theta\subseteq\ker\psi$ it will suffice to show that
$\Theta'\subseteq\ker\psi$, where $\Theta'$ is the set of all words
$$ a=a_1a_2\cdots a_k\in
\Lambdao\Bigl(
\Lambdao\bigl((\pi_\iota(1\otimes M_n(\Cpx))\oup)_{\iota\in I}\bigr),\,
\bigl(\pi_\iota(A_\iota\oup\otimes e)\bigr)_{\iota\in I}\Bigr) $$
such that if $2\le j\le n-1$, if $a_j$ comes from
$$ \Lambdao\bigl((\pi_\iota(1\otimes M_n(\Cpx))\oup)_{\iota\in I}\bigr) $$
and if both $a_{j-1}$
and $a_{j+1}$ come from $\pi_{\iota_1}(A_{\iota_1}\oup\otimes e)$ for some
$\iota_1\in I$, then $a_j\not\in\pi_{\iota_1}(1\otimes M_n(\Cpx))\oup$.
But using that
$$ (1\otimes M_n(\Cpx))(A_\iota\oup\otimes1)
\subseteq(A_\iota\otimes M_n(\Cpx))\oup, $$
we see that every element of $\Theta'$ is equal to an element of
$$ \Lambdao\Bigl(\bigl(\pi_\iota(A_\iota\otimes M_n(\Cpx))\oup)_{\iota\in I}
\Bigr), \tag{\endLo} $$
and the set~(\endLo) of reduced words is a subset of $\ker\psi$ by freeness.
\QED

\proclaim{Lemma \Comprirrat}
Let $I$ be a finite or countably infinite set having at least two elements and
for every $\iota\in I$ let $A(\iota)$ be a type II$_1$ factor.
Let
$$ \MvN=\freeprodi A(\iota) $$
be the free product factor.
Let $0<t<1$ be such that $1/t$ is not an integer
and let $1/t=n+r$ where $n\in\Nats$ and $0<r<1$.
Then $\MvN_t$ is isomorphic to the free product of $|I|+1$ von Neumann algebras
$A(\iota)_t$ ($\iota\in I$) and $\PvN$, with amalgamation over a common two
dimensional subalgebra, $D$, (with respect to trace--preserving conditional
expectations).
These are given by
$$ D=\smd\Cpx r\oplus\smd\Cpx{1-r} $$
and
$$ \PvN=\cases
L(F_{(|I|-1)(t^{-2}-1)+2|I|r(1-r)})&\text{ if }t\le1-\frac1{2|I|} \\
 \vspace{1ex} 
L(F_{2-(|I|+1)(2|I|-1)^{-2}})\oplus\smd\Cpx{2|I|-(2|I|-1)t^{-1}}
 &\text{ if }t>1-\frac1{2|I|},
\endcases \tag{\PvNdetails} $$
where in the case $t>1-\frac1{2|I|}$, the algebra $D$ is embedded into $\PvN$
in such a way that the
projection $1\oplus0\in D$ lies in
$L(F_{2-(|I|+1)(2|I|-1)^{-2}})\oplus0\subseteq\PvN$ and the
projection $0\oplus1\in D$ covers but does not equal the projection
$0\oplus1\in\PvN$.
Equation~(\PvNdetails) should be interpreted to mean $\PvN=L(F_\infty)$ if $I$
is infinite.
\endproclaim
\demo{Proof}
Denote by $\tau$ the tracial state on $\MvN$.
For every $\iota\in I$ let $p(\iota)\in A(\iota)$ be a projection with
$\tau(p(\iota))=t$;
let $v_0(\iota),v_2(\iota),v_3(\iota),\ldots,v_n(\iota)\in A$ be
partial isometries such that
$$ \gather
v_0(\iota)v_0(\iota)^*\le p(\iota), \\
\forall j\in\{2,3,\ldots,n\}\quad v_j(\iota)v_j(\iota)^*=p(\iota) \\
v_0(\iota)^*v_0(\iota)+p(\iota)+\sum_{j=2}^nv_j(\iota)^*v_j(\iota)=1.
\endgather $$
Thus $\tau(v(\iota)_0v(\iota)_0^*)=rt$.
Consider the finite dimensional $*$--subalgebras
$$ B(\iota)=\{v_0(\iota),p,v_2(\iota),v_3(\iota),\ldots,v_n(\iota)\}''\subseteq
A(\iota) $$
and let $\NvN=(\bigcup_{\iota\in I}B(\iota))''$ be the von Neumann algebra they
generate.
Then
$$ B(\iota)\cong M_{n+1}(\Cpx)\oplus M_n(\Cpx) $$
and we choose matrix units
$\bigl(e_{ij}^{(\iota)}\bigr)_{0\le i,j\le n}$ for
$M_{n+1}(\Cpx)\oplus0\subseteq B(\iota)$
and $\bigl(f_{ij}^{(\iota)}\bigr)_{1\le i,j\le n}$ for
$0\oplus M_n(\Cpx)\subseteq B(\iota)$ so that
$$ \align
p(\iota)&=e_{11}^{(\iota)}\oplus f_{11}^{(\iota)},\\
v_0(\iota)&=e_{10}^{(\iota)}\oplus0,\\
\forall j\in\{2,3,\ldots,n\}\quad
 v_j(\iota)&=e_{1j}^{(\iota)}\oplus f_{1j}^{(\iota)}.
\endalign $$

Fix $\iota_0\in I$ and for every $\iota\in I$ write
$p_0(\iota)=e_{11}^{(\iota)}$ and $p_1(\iota)=f_{11}^{(\iota)}$.
Note that $\tau(p_0(\iota))=rt$ while
$\tau(p_1(\iota))=(1-r)t$.
Given $\iota\in I\backslash\{\iota_0\}$ and $j\in\{0,1\}$ and using that
$p_j(\iota)$ and $p_j(\iota_0)$ are free,
that $\tau(p_j(\iota))=\tau(p_j(\iota_0))$ and $\tau$ is
faithful, by~\scite{\DykemaZZFreeDim}{1.1} there is a partial isometry
$u_j(\iota)$ in 
the von Neumann algebra generated by $\{1,p_j(\iota),p_j(\iota_0)\}$ such
that $u_j(\iota)^*u_j(\iota)=p_j(\iota_0)$ and
$u_j(\iota)u_j(\iota)^*=p_j(\iota)$, ($j=0,1$).
Letting $u(\iota)=u_0(\iota)+u_1(\iota)$, we have
$u(\iota)^*u(\iota)=p(\iota_0)$ and $u(\iota)u(\iota)^*=p(\iota)$.
Let also $u(\iota_0)=p(\iota_0)$.
Then
$$ \MvN=\Bigl(\NvN\cup\bigcup_{\iota\in I}p(\iota)A(\iota)p(\iota)\Bigr)''
=\Bigl(\NvN\cup\bigcup_{\iota\in I}u(\iota)^*A(\iota)u(\iota)\Bigr)'', $$
and hence
$$p(\iota_0)\MvN p(\iota_0)=\Bigl(p(\iota_0)\NvN p(\iota_0)\cup
\bigcup_{\iota\in I}u(\iota)^*A(\iota)u(\iota)\Bigr)''. $$
Consider the two dimensional subalgebras
$$ D(\iota)=\Cpx p_0(\iota)+\Cpx p_1(\iota)
\subseteq p(\iota)\NvN p(\iota) $$
and note that $p(\iota)Ap(\iota)\cap p(\iota)\NvN p(\iota)=D(\iota)$.
Henceforth, we usually write simply $D$ instead of $D(\iota_0)$ and $p$ instead
of $p(\iota_0)$.
Note that $u(\iota)^*D(\iota)u(\iota)=D$.
There is a unique trace--preserving conditional expectation
$\Phi:\MvN\to D$
and we will show that with respect to $\Phi$, the
family of $|I|+1$ algebras
$$ p\NvN p,\,
\bigl(u(\iota)^*A(\iota)u(\iota)\bigr)_{\iota\in I} \tag{\freefam} $$
is free with amalgamation over $D$.

Showing that the family~(\freefam)
is free over $D$ amounts to showing that
$$ \Lambdao\biggl(p\NvN p\ominus D,\,
\Bigl(u(\iota)^*A(\iota)u(\iota)\ominus D\Bigr)_{\iota\in I}\biggr)
\subseteq\ker\Phi, $$
which is easily seen to be equivalent to
$$ \Lambdao\biggl(p\NvN p\ominus D,\,
\Bigl(u(\iota)^*\Bigl(p(\iota)A(\iota)p(\iota)\ominus D(\iota)\Bigr)
u(\iota)\Bigr)_{\iota\in I}\biggr)
\subseteq\ker\tau. \tag{\inkertau} $$
Let $x$ belong to the left--hand--side of~(\inkertau).
Note that
$$ p(\iota)A(\iota)p(\iota)\ominus D(\iota)=
 \bigl(p_0(\iota)A(\iota)p_0(\iota)\bigr)\oup
 +p_0(\iota)A(\iota)p_1(\iota)+p_1(\iota)A(\iota)p_0(\iota)
 +\bigl(p_1(\iota)A(\iota)p_1(\iota)\bigr)\oup $$
and that $u(\iota)\in\NvN$ and
$u(\iota)(p\NvN p\ominus D)u(\iota)^*=p(\iota)\NvN p(\iota)\ominus D(\iota)$.
Let
$$ \Omegao_\iota\,=\,\bigl(p_0(\iota)A(\iota)p_0(\iota)\bigr)\oup
\;\;\cup\;\;p_0(\iota)A(\iota)p_1(\iota)
\;\;\cup\;\;p_1(\iota)A(\iota) p_0(\iota)
\;\;\cup\;\;\bigl(p_1(\iota)A(\iota)p_1(\iota)\bigr)\oup.  $$
Separating off the $u(\iota)^*$ and $u(\iota)$ appearing in $x$ and multiplying
some neighbors,
we see that if $x\not\in p\NvN p\ominus D$ then $x=x'$ for some $x'\in\Theta$,
where $\Theta$ is the set of all words, $w=w_1w_2\cdots w_m$, belonging to
$\Lambdao\bigl(\NvN,(\Omegao_\iota)_{\iota\in I}\bigr)$
such that
\roster
\item"$\bullet$"
 $w$ has at least one letter from some $\Omegao_\iota$;
\item"$\bullet$"
 if a letter, $w_j$ in $w$ comes from $\NvN$ and if it lies between two
 letters, $w_{j-1}$ and $w_{j+1}$, both of which belong to the same
 $\Omegao_\iota$ then $w_j\in p(\iota)\NvN p(\iota)\ominus D(\iota)$.
\endroster
Now
$$ \NvN=\clspan\Bigl(\{1\}\cup
 \Lambdao\bigl((B(\iota)\oup)_{\iota\in I}\bigr)\Bigr), $$
where $\clspan$ means the closure of the linear span in the w$^*$--topology.
By Kaplansky's density theorem, it follows that every element of $\NvN$ is
the limit in strong operator topology of a bounded net in
$\lspan\bigl(\{1\}\cup\Lambdao((B(\iota)\oup)_{\iota\in I})\bigr)$.
Furthermore, since $p(\iota)B(\iota)p(\iota)\ominus D(\iota)=\{0\}$, it
follows, fixing $\iota'\in I$,
that every element of $p(\iota')\NvN p(\iota')\ominus D(\iota')$ is the
s.o.--limit of a bounded net in
$$ \lspan\Bigl(\Lambdao\bigl((B(\iota)\oup)_{\iota\in I}\bigr)
 \bigm\backslash B(\iota')\oup\Bigr). $$
Hence it will suffice to show that $\Theta'\subseteq\ker\tau$ where $\Theta'$ is
the set of all words, $w=w_1w_2\cdots w_m$ belonging to
$$ \Lambdao\Bigl(\Lambdao\bigl((B(\iota)\oup)_{\iota\in I}\bigr),
 (\Omegao_\iota)_{\iota\in I}\Bigr) $$
such that
\roster
\item"$\bullet$"
 $w$ has at least one letter from some $\Omegao_\iota$;
\item"$\bullet$"
 if a letter $w_j$ in $w$ comes from
 $\Lambdao\bigl((B(\iota)\oup)_{\iota\in I}\bigr)$ and if it lies
 between two letters, $w_{j-1}$ and $w_{j+1}$, both of which belong to
 $\Omegao_{\iota'}$ for the same fixed $\iota'\in I$ then
 $w_j\notin B(\iota')\oup$.
\endroster
Let $y\in\Theta'$.
Since $p(\iota)=p_0(\iota)+p_1(\iota)$, since $p_0(\iota)$ and $p_1(\iota)$ are
minimal projections in $B(\iota)$ and since 
$p_0(\iota)B(\iota)p_1(\iota)=\{0\}$ it follows that
$B(\iota)\Omegao_\iota B(\iota)\subseteq A(\iota)\oup$.
Hence, writing out all letters in $y$ that come from
$\Lambdao\bigl((B(\iota)\oup)_{\iota\in I}\bigr)$ and
combining some neighbors, we see that $y=y'$ where $y'$ is a word belonging to
$\Lambdao\bigl((A(\iota)\oup)_{\iota\in I}\bigr)$.
But now by freeness of $(A(\iota))_{\iota\in I}$ we have $\tau(y')=0$.
This finishes the proof that the family~(\freefam) is free over $D$.

Let $\PvN=p\NvN p$.
It remains to describe $D$, $\PvN$ and the embedding of $D$ in $\PvN$ in more
detail.
In the notation of~\cite{\DykemaZZFreeDim}, we have
$$ B(\iota)=\smdp{M_{n+1}(\Cpx)}{(n+1)rt}{e_{ij}^{(0)}}
 \oplus\smdp{M_n(\Cpx)}{n(1-r)t}{e_{ij}^{(1)}}, \tag{\Bpic}$$
and $\NvN=\freeprodi B(\iota)$, taking the free product with respect
to the traces indicated in~(\Bpic).
This free product is easily computed using the results
of~\cite{\DykemaZZFreeDim}.
We find that
\vskip2ex
$$ \NvN=\cases
L(F_{|I|(1-t^2+2t^2r(1-r))})&
 \text{ if }n\ge2
 \text{ or }\Bigl(n=1\text{ and }t\le1-\frac1{2|I|}\Bigr) \\ \vspace{2ex}
L(F_{2-\frac5{4|I|}})\oplus\smdp\Cpx{1-2|I|(1-t)}q
 &\text{ if }n=1\text{ and }t>1-\frac1{2|I|},
\endcases $$
where $q=\bigwedge_{\iota\in I}p_1(\iota)$.
Since $p=e_{11}^{(\iota)}\oplus f_{11}^{(\iota)}$ and has trace $t$,
compressing we get~(\PvNdetails).
Now it is clear that $D$ embeds into $\PvN$ as described.
\QED

\proclaim{Lemma \Stickout}
Let $A$ be a II$_1$--factor and denote its tracial state by $\tau_A$.
Let $B$ be a von Neumann algebra with faithful tracial state $\tau_B$, and
suppose that $p$ is a minimal and central projection in $B$ with
$\tau_B(p)=1/n$ for some integer $n\ge2$.
Let
$$ (\MvN,\tau)=(A,\tau_A)*(B,\tau_B). $$
Then $\MvN$ is a II$_1$--factor.
Let $\tau^{(1/n)}$ (respectively $\tau_A^{(1/n)}$) denote the tracial state on
$\MvN_{1/n}$, (respectively $A_{1/n}$).
Then
$$ (\MvN_{1/n},\tau^{(1/n)})\cong(A_{1/n},\tau_A^{(1/n)})*(pBp,n\tau_B\restrict_{pBp})
*(L(F_{2(n-1)}),\sigma), \tag{\StickoutEq} $$
where $\sigma$ is the tracial state on the indicated free group factor.
\endproclaim
\demo{Proof}
Resorting to the notation of~\cite{\DykemaZZFreeDim}, we have
$$ \align
A&=A_{1/n}\otimes M_n(\Cpx) \\
B&=pBp\oplus\smd{\Cpx}{(n-1)/n},
\endalign $$
and we consider the von Neumann subalgebras
$$ \matrix
\MvN&=&\bigl(&A_{1/n}&\otimes&M_n(\Cpx)&\bigr)&
*&\bigl(&pBp&\oplus&\smdp{\Cpx}{(n-1)/n}{1-p}&\bigr) \\\vspace{1ex}
\cup \\\vspace{1ex}
\NvN&=&\bigl(&A_{1/n}&\otimes&M_n(\Cpx)&\bigr)&
*&\bigl(&\smdp\Cpx{}p&\oplus&\smdp{\Cpx}{(n-1)/n}{1-p}&\bigr) \\\vspace{1ex}
\cup \\\vspace{1ex}
\PvN&=&\bigl(&\Cpx&\otimes&M_n(\Cpx)&\bigr)&
*&\bigl(&\smdp\Cpx{}p&\oplus&\smdp{\Cpx}{(n-1)/n}{1-p}&\bigr).
\endmatrix $$
By~\cite{\DykemaZZFreeDim}, $\PvN\cong L(F_{(n^2+2n-3)/n^2})$ and hence
$\PvN_{1/n}\cong L(F_{2(n-1)})$;
a double application of Theorem~1.2 of~\cite{\DykemaZZFreeDim} shows that
$p\MvN p$ is generated by free copies of $A_{1/n}$, $\PvN_{1/n}$ and $pBp$.
Thus $p\MvN p$ is a factor, as is $\MvN$, and~(\StickoutEq) holds.
\QED

\proclaim{Lemma \OnlyFreeDim}
Let $\NvN$ be a II$_1$--factor and let $\tau_\NvN$ be its tracial state.
Let $A_i$ be a finite von Neumann algebra that is either finite dimensional,
approximately finite dimensional, an interpolated free group factor or a direct
sum of these, and let $\tau_{A_i}$ be a faithful tracial state on $A_i$,
($i=1,2$).
Let
$$ (\MvN_i,\tau_i)=(\NvN,\tau_\NvN)*(A_i,\tau_{A_i}). $$
If
$$ \fdim(A_1,\tau_{A_1})=\fdim(A_2,\tau_{A_2}) \tag{\eqfdim} $$
then $\MvN_1\cong\MvN_2$ and $\MvN_1$ is a factor.
\endproclaim
\demo{Proof}
If $\fdim(A_i,\tau_{A_i})=0$ then $A_i=\Cpx$ and $\MvN_i=\NvN$.
Suppose $\fdim(A_i,\tau_{A_i})>0$.
Then each $A_i$ has linear dimension at least two;
let
$$ (\PvN_i^{(n)},\tau_{\PvN_i}^{(n)})=(M_n(\Cpx),\tr_n)*(A_i,\tau_{A_i}), $$
($i=1,2$);
then by~\cite{\DykemaZZFreeDim}, for $n\in\Nats$ large enough both
$\PvN_1^{(n)}$ and $\PvN_2^{(n)}$ are factors;
moreover, under the assumption~(\eqfdim), $\PvN_1\cong\PvN_2$.
By Theorem~1.2 of~\cite{\DykemaZZFreeDim},
$(\MvN_i)_{1/n}\cong\NvN_{1/n}*(\PvN_i)_{1/n}$;
hence $\MvN_1\cong\MvN_2$ is a factor.
\QED

\proclaim{Theorem \Compr}
Let $I$ be a finite or countably infinite set having at least two elements and
for each $\iota\in I$ let $A(\iota)$ be a type II$_1$ factor.
Consider the free product factor $\MvN=\freeprodi A(\iota)$
and let $t$ be a real number satisfying $0<t<1$.
Then
$$
\MvN_t\cong\cases
\bigl(\freeprodi A(\iota)_t\bigr)*L(F_{(|I|-1)(t^{-2}-1)})&
 \text{if }(|I|-1)(t^{-2}-1)>1, \\ \vspace{1ex} 
\bigl(\freeprodi A(\iota)_t\bigr)*R&
 \text{if }(|I|-1)(t^{-2}-1)=1, \\ \vspace{1ex}
\bigl(\freeprodi A(\iota)_t\bigr)*(R\oplus\smd\Cpx\alpha)&
 \text{if }(|I|-1)(t^{-2}-1)<1,
\endcases
\tag{\cutIfinite} $$
where $R$ is the approximately finite dimensional (i.e. hyperfinite)
II$_1$--factor, in the last instance $1-\alpha^2=(|I|-1)(t^{-2}-1)$ and
the notation $R\oplus\smd\Cpx\alpha$ means that we are taking the free product
with respect to the trace on $R\oplus\Cpx$ that takes the value $\alpha$ on
the element $0\oplus1$.

If $I$ is infinite then
$$ \align
\MvN&\cong\MvN*L(F_\infty) \tag{\stable} \\
\MvN_t&\cong\freeprodi A(\iota)_t. \tag{\cutIinfinite}
\endalign $$
\endproclaim
\demo{Proof}
We first prove~(\cutIfinite).
If $t=1/n$ for some integer $n\ge2$ then we are in the first instance
of~(\cutIfinite) and the isomorphism follows directly from Lemma~\Comproneonn.
Assume now that $t$ is not a reciprocal integer.
We first prove the case when $(|I|-1)(t^{-2}-1)>1$.
{}From Lemma~\Comprirrat{} we have that $\MvN_t$ is isomorphic to the free
product of the family of $|I|+1$ von Neumann algebras
$$ L(F_x),\quad\bigl(A(\iota)_t\bigr)_{\iota\in I} $$
all amalgamated over copies of
$$ D=\smd\Cpx{1-r}\oplus\smd\Cpx r, $$
where $x=(|I|-1)(t^{-2}-1)+2|I|r(1-r)$.
Note that the free dimension of $D$ is $2r(1-r)$.
But by~\cite{\DykemaZZFreeDim}, $L(F_x)$ is isomorphic to the
free product of $L(F_y)$ and $|I|$ copies of $D$, where $y=(|I|-1)(t^{-2}-1)$.
Letting these copies of $D$ be denoted $(D_\iota)_{\iota\in I}$ we thus have
$$ L(F_x)\cong L(F_y)*(\freeprodi D_\iota). $$
Let $B(\iota)$ be the von Neumann algebra generated by $A(\iota)_t$ and
$D_\iota$, which are free.
Then $\MvN_t$ is isomorphic to the free product
of $(|I|+1)$ II$_1$--factors $L(F_y)$ and $\bigl(B(\iota)\bigr)_{\iota\in I}$,
all amalgamated over embedded copies $D\hookrightarrow L(F_y)$ and
$D\hookrightarrow B(\iota)$.
Ostensibly, our embedded copy of $D$ in $B(\iota)$ lies inside of
$A(\iota)_t\subseteq B(\iota)$.
But $B(\iota)$ is a factor by Lemma~\OnlyFreeDim, and since $D$ is finite
dimensional and abelian, it follows that every two (tracially identical) copies
of $D$ inside of $B(\iota)$ are equivalent via an inner automorphism.
Hence we know that our embedded copy of $D\hookrightarrow B(\iota)$ is freely
complemented in $B(\iota)$ by a subalgebra isomorphic to $A(\iota)_t$.
Therefore $\MvN_t\cong L(F_y)*\bigl(\freeprodi A(\iota)_t\bigr)$, as required.

If $(|I|-1)(t^{-2}-1)=1$ then the proof goes exactly as above but where we
replace $L(F_y)$ everywhere by $R$.

Now suppose that $y\eqdef(|I|-1)(t^{-2}-1)<1$.
By choice of $\alpha$, the free dimension of $R\oplus\smd\Cpx\alpha$ is $y$.
Let $\beta>\alpha$ be such that $1-\beta=1/m$ for some positive integer $m$.
When $z>1$ the free dimension of $L(F_z)\oplus\smd\Cpx\beta$ is
$1+(1-\beta)^2(z-1)-\beta^2$, and we choose $z=m^2y-2m+2$ so that the free
dimension of $L(F_z)\oplus\smd\Cpx\beta$ is equal to $y$.
By Lemma~\OnlyFreeDim, we have
$$ \bigl(\freeprodi A(\iota)_t\bigr)*(R\oplus\smd\Cpx\alpha)
\cong\bigl(\freeprodi A(\iota)_t\bigr)*(L(F_z)\oplus\smd\Cpx\beta), $$
and we will now show that $\MvN_t$ is isomorphic to the right--hand--side, by
rescaling both by $1/m$.
Since $t/m<1/m\le1/2$, when rescaling $\MvN$ by $t/m$ we may use the first
instance of~(\cutIfinite) to get
$$ (\MvN_t)_{1/m}=\MvN_{t/m}
\cong\bigl(\freeprodi A(\iota)_{t/m}\bigr)*L(F_{(|I|-1)(m^2t^{-2}-1)}). $$
On the other hand, using Lemma~\Stickout{} and then the first instance
of~(\cutIfinite) we find that
$$ \align
\biggl(\bigl(\freeprodi A(\iota)_t\bigr)
 *(L(F_z)\oplus\smd\Cpx\beta)\biggr)_{1/m}
&\cong\bigl(\freeprodi A(\iota)_t\bigr)_{1/m}*L(F_z)*L(F_{2(m-1)}) \\
&\cong\bigl(\freeprodi A(\iota)_{t/m}\bigr)*L(F_{(|I|-1)(m^2-1)})
 *L(F_{z+2(m-1)}) \\
&\cong\bigl(\freeprodi A(\iota)_{t/m}\bigr)*L(F_{(|I|-1)(m^2t^{-2}-1)}).
\endalign $$
This completes the proof of~(\cutIfinite).

In order to prove~(\stable), we will show that
$$ \MvN_{1/2}\cong\Bigl(\MvN*L(F_\infty)\Bigr)_{1/2}. $$
But using~(\cutIfinite) we see that both sides are isomorphic to
$$ \bigl(\freeprodi A(\iota)_{1/2}\bigr)*L(F_\infty). $$
Now~(\cutIinfinite) follows from~(\cutIfinite) and~(\stable).
\QED

Recall that $\Feu(\MvN)$ denotes the fundamental group of a type II$_1$ factor
$\MvN$.
{}From~(\cutIinfinite), information about the fundamental group of the free
product of infinitely many type II$_1$ factors can now be obtained.
\proclaim{Corollary \FGpinfinite}
For every $n\in\Nats$ let $A(n)$ be a type II$_1$ factor and let
$$ \MvN=\freeprod{n=1}\infty A(n) $$
be their free product factor.
Then
$$ \bigcap_{n=1}^\infty\Feu(A(n))\subseteq\Feu(\MvN). $$
\endproclaim

\heading
\S\secAutom.  Semifinite free products and trace--scaling automorphisms
\endheading

The following construction is an example of a realization of a free product
situation with a von Neumann algebra $\MvN$ endowed with a semifinite
weight $\phi$.
We only know how to do this
  \footnote{After this was written, two persons, Roland Speicher and a referee,
  have informed us that it is possible to perform this construction letting $X$
  be, for example, a self--adjoint element having infinitely divisible
  distribution.}
when we let the other
algebra be generated by a fixed semicircular element, $X$.
The free product weight is not finite if the original weight
$\phi$ is not finite, and then
the fixed semicircular element is not an element of the free product but rather
gives rise to an unbounded quadratic form.
Such a method was first used in~\cite{\RadulescuZZOneParamGp}, on which the
proofs here are based.

\proclaim{Theorem and Definition \fpSX}
Let $\MvN$ be a von Neumann algebra and let $\phi$ be a normal, faithful,
semifinite weight on $\MvN$.
Let $\Pc_0$ be the set of all projections $p\in\MvN$ such that
$\phi(p)<\infty$.
Let $X$ be an indeterminate and let $\Cpx[X]$ be the algebra of polynomials in
$X$ endowed with the $*$--operation so that $X^*=X$.
Consider the unital $*$--algebra free product $\Afr=\MvN*\Cpx[X]$ and for every
$p\in\Pc_0$ consider the $*$--subalgebra $\Ac_p\subseteq\Afr$ defined by
$$ \Ac_p=p\MvN p+\lspan\{a_0Xa_1Xa_2\cdots Xa_n\mid n\ge1,\,a_j\in p\MvN p\}; $$
let $\Ac=\bigcup_{p\in\Pc_0}\Ac_p$.
Then there is a unique positive linear functional $\psi:\Ac\to\Cpx$ such that
for every $p\in\Pc_0$,
\roster
\item"(i)" $\psi\restrict_{p\MvN p}=\phi\restrict_{p\MvN p}$;
\item"(ii)" with respect to $\phi(p)^{-1}\psi\restrict_{\Ac_p}$, $pXp$ is a
 semicircular element with second moment $\phi(p)$;
\item"(iii)" $\{pXp\}$ and $p\MvN p$ are free with respect to
 $\phi(p)^{-1}\psi\restrict_{\Ac_p}$.
\endroster
Moreover, $\Ac$ acts by bounded operators on $L^2(\Ac,\psi)$ via the GNS
representation, $\pi_\psi$.
Let $\NvN=\pi_\psi(\Ac)''$ be the von Neumann algebra generated by the image of
$\pi_\psi$.
Then $\psi$ extends to a normal, faithful, semifinite weight on $\NvN$, also
denoted $\psi$, and there is a canonical normal injective $*$--homomorphism
$\MvN\hookrightarrow\NvN$ such that $\psi\restrict_\MvN=\phi$.
If $\phi$ is a trace then $\psi$ is a trace.

We denote this construction by
$$ (\NvN,\psi)=(\MvN,\phi)*\Seu X. $$
\endproclaim
\demo{Proof}
Consider the state $\phit_p$ on $\MvN$ given by
$\phit_p(a)=\phi(p)^{-1}\phi(pap)$;
let $\sigma_p$ be the state on the $*$--algebra $\Cpx[X]$ so that $X$ becomes a
semicircular element with second moment $\phi(p)$, i.e\.
$$ \sigma_p(X^k)=\frac1{2\pi\phi(p)}\int_{-2\sqrt{\phi(p)}}^{2\sqrt{\phi(p)}}
 t^k\sqrt{4\phi(p)-t^2}\dif t. $$
Let $\psit_p=\phit_p*\sigma_p$ be the resulting free product state on $\Afr$
and let $\psi_p=\phi(p)\psit_p$.

We will show that if $p,q\in\Pc_0$ then $\psi_p$ and $\psi_q$ agree on
$\Ac_p\cap\Ac_q$.
To show this we may without loss of generality assume $p\le q$, which implies
$\Ac_p\subseteq\Ac_q$.
It is clear that
$\psi_p\restrict_{p\MvN p}=\phi\restrict_{p\MvN p}=\psi_q\restrict_{p\MvN p}$,
and we must show that $\psi_p$ and $\psi_q$ agree on elements of the form
$a_0Xa_1Xa_2\cdots Xa_n$ for $n\ge1$ and $a_j\in p\MvN p$.
For this we will use R\. Speicher's combinatorial treatment of
freeness~\cite{\SpeicherZZNoncrossing}; we
will use a slight modification of Theorem~2.17{} of~\cite{\SpeicherZZFreeCal},
(alternatively, see~\S3.4 of~\cite{\SpeicherZZComb}),
which implies that
$$ \psit_p(a_0Xa_1\cdots Xa_n)=
\sum_{\pi\in\NC(\{2,4,\ldots,2n\})}
k_{\pi}[\undersetbrace{n\text{ times}}\to{X,X,\ldots,X}]
(\phit_p)_{\pi^c}[a_0,a_1,\ldots,a_n]. $$
Here $\NC(\cdot)$ is the set of all non--crossing partitions on an ordered set;
$k_\pi$ is the cummulant function associated to the state $\sigma_p$ on
$\Cpx[X]$;
$\pi^c\in\NC(\{1,3,5,\ldots,2n+1\})$ is the Kreweras
complement~\cite{\Kreweras} of the non--crossing partition $\pi$, namely, it is
the largest non--crossing partition of $\{1,3,5,\ldots,2n+1\}$ such that
$\pi\cup\pi^c$ is a non--crossing partition of
$\{1,2,\ldots,2n+1\}$;
finally, if $\pi^c=\{B_1,B_2,\ldots,B_k\}$ where
$B_j=\{i^{(j)}_1,i^{(j)}_2,\ldots,i^{(j)}_{\ell(j)}\}$ are the blocks of the
partition with $i^{(j)}_1<i^{(j)}_2<\cdots<i^{(j)}_{\ell(j)}$,
and if $\alpha$ is the order preserving bijection
$\{1,3,\ldots,2n+1\}\to\{0,1,\ldots,n\}$, then
$$ (\phit_p)_{\pi^c}[a_0,a_1,\ldots,a_n]=
\prod_{j=1}^k\phit_p(a_{\alpha(i_1^{(j)})}a_{\alpha(i_2^{(j)})}
\cdots a_{\alpha(i_{\ell(j)}^{(j)})}). \tag{\phitp} $$
But since $X$ is semicircular with respect to $\sigma_p$ with second moment
$\phi(p)$, it follows that
$k_\pi[X,X,\ldots,X]$ is zero unless $n$ is even and $\pi$ is a pairing
i.e\. all blocks of $\pi$ have two elements, and in that case
$k_\pi[X,X,\ldots,X]=\phi(p)^{n/2}$.
However, it is an easy lemma that if $\pi$ is a non--crossing pairing then
$\pi^c$ has exactly $\frac n2+1$ blocks.
Therefore, denoting by $\NCP(\cdot)$ the set of all non--crossing pairings of
an ordered set, and letting $(\phi)_{\pi^c}[a_0,a_1,\ldots,a_n]$ be defined
analogously to~(\phitp),
we have
$$ \aligned
\psi_p(a_0Xa_1\cdots Xa_n)&=\phi(p)\psit_p(a_0Xa_1\cdots Xa_n) \\
&=\phi(p)\sum_{\pi\in\NCP(\{2,4,\ldots,2n\})}\phi(p)^{n/2}
(\phit_p)_{\pi^c}[a_0,a_1,\ldots,a_n] \\
&=\sum_{\pi\in\NCP(\{2,4,\ldots,2n\})}(\phi)_{\pi^c}[a_0,a_1,\ldots,a_n].
\endaligned \tag{\psipexpr} $$
But this last expression does not depend explicitly on $p$, and we get the
same expression for $\psi_q(a_0Xa_1\cdots Xa_n)$.
Hence $\psi_p$ and $\psi_q$ agree on $\Ac_p\cap\Ac_q$.

We may therefore define $\psi:\Ac\to\Cpx$ by $\psi\restrict_{\Ac_p}=\psi_p$.
Since every $\psi_p$ is positive, also $\psi$ is positive.
By construction,~(i), (ii) and~(iii) hold, and the uniqueness of $\psi$ is
clear.
Given $a\in\Ac$, we denote its corresponding element in $L^2(\Ac,\psi)$ by
$\ah$.
The GNS representation of $\Ac$ on $L^2(\Ac,\psi)$ is defined by
$\pi_\psi(b)\ah=(ba)\hat{\;}$.
In order to show that the GNS representation $\pi_\psi(\cdot)$ acts by bounded
operators on $L^2(\Ac,\psi)$, it will suffice to show for arbitrary $q\in\Pc_0$
and  $a_0,a_1\in q\MvN q$ that $\pi_\psi(a_0Xa_1)$ is bounded.
Since $L^2(\Ac,\psi)=\bigcup_{p\in\Pc_0}L^2(\Ac_p,\psi_p)$, it will suffice to
show that $\nm{\pi_{\psi_p}(a_0Xa_1)}$ is uniformly bounded for $p\in\Pc_0$,
$q\le p$, where $\pi_{\psi_p}(a_0Xa_1)$ acts on $\L^2(\Ac_p,\psi_p)$.
But the unitary $W_p:L^2(\Ac_p,\psi_p)\to L^2(\Ac_p,\psit_p)$ given by
$\ah\to\phi(p)^{1/2}\ah$ intertwines $\pi_{\psi_p}$ and $\pi_{\psit_p}$, so we
must only show that $\nm{\pi_{\psit_p}(a_0Xa_1)}$ is uniformly bounded.
Now with respect to $\psit_p$, $X$ is a semicircular element with second moment
$\phi(p)$, and  $q$ is a projection free from $X$ and with
$\psit_p(q)=\phi(q)/\phi(p)$.
Therefore, by~\cite{\VoiculescuZZCircSemiCirc}, $qXq$ is with respect to
the faithful state $\frac{\phi(p)}{\phi(q)}\psit_p\restrict_{q\Ac_pq}$ a semicircular element with second moment
$\phi(q)$, and hence $\nm{\pi_{\psi_p}(qXq)}=2\phi(q)^{1/2}$.
Therefore $\nm{\pi_{\psi_p}(a_0Xa_1)}\le2\phi(q)^{1/2}\nm{a_0}\,\nm{a_1}$ uniformly in
$p$.

Let $\NvN=\pi_\psi(\Ac)''$.
For $p\in\Pc_0$ let $E_p$ be the orthogonal projection of $L^2(\Ac,\psi)$ onto
its subspace $L^2(\Ac_p,\psi_p)$;
note that if $q\in\Pc_0$, $p\le q$ and $a_0,a_1,\ldots,a_n\in q\MvN q$ then
$$ E_p\pi_\psi(a_0Xa_1\cdots Xa_n)\restrict_{L^2(\Ac_p,\psi_p)}
=\pi_{\psi_p}(pa_0pXpa_1pXpa_2p\cdots Xpa_np), $$
so $a\mapsto E_paE_p$ defines a normal, conditional expectation from $\NvN$
onto $\pi_{\psi_p}(\Ac_p)''$.
We can define $\psi$ on $\NvN$ by $\psi(x)=\sup_{p\in\Pc_0}\psi_p(E_pxE_p)$,
which gives a normal, semifinite weight on $\NvN$.
To see that $\psi$ is faithful, it suffices to note that $E_p\nearrow1$
as $p\nearrow1$ and that Voiculescu
proved~\cite{\VoiculescuZZSymmetries} that each $\psi_p$ is faithful on
$\pi_{\psi_p}(\Ac_p)''$.
The embeddings $p\MvN p\hookrightarrow\pi_{\psi_p}(\Ac_p)''$ arising from the
free product construction taken together give an embedding
$\MvN\hookrightarrow\NvN$ such that $\psi$ restricts to $\phi$.
\QED

\proclaim{Theorem \IIinf}
Let $I$ be a finite or countably infinite set and for every $\iota\in I$ let
$A(\iota)$ be a type II$_1$ factor with tracial state $\tau_{A(\iota)}$;
letting $\Hil$ be separable, infinite
dimensional Hilbert space, consider the type II$_\infty$ von Neumann algebra
$$ \MvN=\bigoplus_{\iota\in I}A(\iota)\otimes B(\Hil). $$
Let $\Tr$ be the trace on $B(\Hil)$ taking value $1$ on minimal projections,
and consider the normal, faithful, semifinite, tracial weight
$$ \tau_\MvN=\bigoplus_{\iota\in I}\tau_{A(\iota)}\otimes\Tr $$
on $\MvN$.
Let
$$ (\NvN,\tau)=(\MvN,\tau_\MvN)*\Seu X. $$
Then $\NvN$ is the II$_\infty$ factor isomorphic to $A\otimes B(\Hil)$, where
$A$ arises as the free product of II$_1$--factors
$$ (A,\tau_A)=(L(F_\infty),\tau_{F_\infty})*
\Bigl(\freeprodi(A(\iota),\tau_{A(\iota)})\Bigr). \tag{\Afp} $$
If for every $\iota\in I$, $A(\iota)\otimes B(\Hil)$ has a one--parameter group
of trace scaling automorphisms, then $\NvN$ has a one--parameter group of
trace--scaling automorphisms.
\endproclaim
\demo{Proof}
We abuse notation by equating $A(\iota)\otimes B(\Hil)$ with the corresponding
direct summand subalgebra of $\MvN$.
Let $(e_{jk})_{j,k\in\Nats}$ be a system of matrix units for $B(\Hil)$; for
$\iota\in I$ let $1_{A(\iota)}$ be the identity element of $A(\iota)$ and let
$e_{jk}^{(\iota)}=1_{a(\iota)}\otimes e_{jk}$.
Choose $\iota_0\in I$, and for $\iota\in I\backslash\{\iota_0\}$ let $u_\iota$
be the polar part of
$e_{11}^{(\iota)}Xe_{11}^{(\iota_0)}$.
Letting $p=e_{11}^{(\iota)}\oplus e_{11}^{(\iota_0)}$,
recalling that $X$ is with respect to $\psit_p$ a semicircular element free
from $e_{11}^{(\iota_0)}$, noting that $\psit_p(e_{11}^{(\iota_0)})=1/2$ and
using~\cite{\VoiculescuZZCircSemiCirc}, we see that $u_\iota$ is a partial
isometry satisfying $u_\iota^*u_\iota=e^{(\iota_0)}_{11}$ and
$u_\iota u_\iota^*=e^{(\iota)}_{11}$.
Let also $u_{\iota_0}=e^{(\iota_0)}_{11}$ and for $\iota_1,\iota_2\in I$ and
$j,k\in\Nats$ let
$$ f(\iota_1,j;\iota_2,k)
=e^{(\iota_1)}_{j,1}u_{\iota_1}u_{\iota_2}^*e^{(\iota_2)}_{1,k}. $$
Then
$$ f(\iota_1,j;\iota_2,k)f(\iota_1',j';\iota_2',k')
=\delta_{\iota_2,\iota_1'}\delta_{k,j'}f(\iota_1,j;\iota_2',k'), $$
which shows that
$\bigl(f(\iota_1,j;\iota_2,k)\bigr)_{\iota_1,\iota_2\in I,\,j,k\in\Nats}$ is a
system of matrix units whose supports sum to the identity element of $\NvN$.
Note that $\tau\bigl(f(\iota_0,1;\iota_0,1)\bigr)=1$.
Let $\QvN=f(\iota_0,1;\iota_0,1)\NvN f(\iota_0,1;\iota_0,1)$.
We will show that $\QvN$ is isomorphic to the II$_1$--factor $A$, arising as the
free product~(\Afp).
Let $D(\iota)=A(\iota)\otimes e_{11}$.
Now $\NvN$ is generated by the set of matrix units
$$ \{f(\iota_1,j;\iota_2,k)\mid\iota_1,\iota_2\in I,\,j,k\in\Nats\}
\tag{\matrixunits} $$
together with
$$ \{e_{jj}^{(\iota_1)}Xe_{kk}^{(\iota_2)}\mid
\iota_1,\iota_2\in I,\,j,k\in\Nats\}\cup
\bigcup_{\iota\in I}D(\iota). \tag{\generators} $$
Let $>$ be a well-ordering of $I$ whose minimal element is $\iota_0$ and denote
also by $>$ the resulting lexicographic ordering of $I\times\Nats$.
Using the matrix units~(\matrixunits) to pull the generators~(\generators) back
to $f(\iota_0,1;\iota_0,1)$, we see that $\QvN$ is generated by
$$ \align
&\{f(\iota_0,1;\iota,j)Xf(\iota,j;\iota_0,1)\mid\iota\in I,\,j\in\Nats\}\cup \\
&\cup\bigl\{f(\iota_0,1;\iota,1)Xf(\iota_0,1;\iota_0,1)
 \mid\iota\in I\backslash\{\iota_0\}\bigr\}\cup \\
&\cup\{f(\iota_0,1;\iota,j)Xf(\iota_0,1;\iota_0,1)
 \mid\iota\in I,\,j>1\bigr\}\cup \\
&\cup\bigl\{f(\iota_0,1;\iota_1,j)Xf(\iota_2,k;\iota_0,1)
 \mid\iota_1,\iota_2\in I,\,j,k\in\Nats,\,
 (\iota_1,j)>(\iota_2,k)>(\iota_0,1)\bigr\}\cup \\
&\cup\bigcup_{\iota\in I}u_\iota^*D(\iota)u_\iota.
\endalign $$
We see that for $\iota\in I\backslash\{\iota_0\}$,
$f(\iota_0,1;\iota,1)Xf(\iota_0,1;\iota_0,1)=b_i$, where $b_\iota$
is the positive part of $e_{11}^{(\iota)}Xe_{11}^{(\iota_0)}$.
We will show that with respect to the state $\psi\restrict_\QvN$,
$$  F_1\eqdef\bigl(f(\iota_0,1;\iota,j)Xf(\iota,j;\iota_0,1)
 \bigr)_{\iota\in I,\,j\in\Nats} $$
is a free family of semicircular elements,
$$  F_2\eqdef(b_\iota)_{\iota\in I\backslash\{\iota_0\}} $$
is a free family of quartercircular elements,
$$  F_3\eqdef
\bigl(f(\iota_0,1;\iota,j)Xf(\iota_0,1;\iota_0,1) \bigr)_{\iota\in I,\,j>1} $$
is a $*$--free family of circular elements,
$$  F_4\eqdef\bigl(f(\iota_0,1;\iota_1,j)Xf(\iota_2,k;\iota_0,1)
 \bigr)_{\iota_1,\iota_2\in I,\,
 j,k\in\Nats,\,(\iota_1,j)>(\iota_2,k)>(\iota_0,1)} $$
is a $*$--free family of circular elements,
$$  F_5\eqdef\bigl(u_\iota^*D(\iota)u_\iota\bigr)_{\iota\in I} $$
is a free family of subalgebras of $\QvN$ and all the families taken together,
$$ ( F_1, F_2, F_3, F_4, F_5), $$
form a $*$--free family.
This will suffice to prove that $\QvN$ is isomorphic to $A$, because the
families $ F_1$, $ F_2$, $ F_3$ and $ F_4$ taken together generate a
copy of $L(F_\infty)$, while the family $ F_5$ consists of one copy of each
$A(\iota)$.

Given a finite subset $E\subseteq I$ such that
$\iota_0\in E$, and given  $K\in\Nats$, let $ F_j^{E,K}$ be defined like
$ F_j$, but quantifying over $E$ instead of $I$ and over $\{1,2,\ldots,K\}$
instead of $\Nats$.
It will suffice to show that
\roster
\item"(i)" $ F_1^{E,K}$ is a free family of semicircular elements;
\item"(ii)" $ F_2^{E,K}$ is a free family of quartercircular elements;
\item"(iii)" $ F_3^{E,K}$ is a $*$--free family of circular elements;
\item"(iv)" $ F_4^{E,K}$ is a $*$--free family of circular elements;
\item"(v)" $ F_5^{E,K}$ is a free family of subalgebras of $\QvN$;
\item"(vi)" the family
 $( F_1^{E,K}, F_2^{E,K}, F_3^{E,K}, F_4^{E,K}, F_5^{E,K})$ is
 $*$--free.
\endroster
Consider
$$ p\eqdef\sum_{\iota\in E}\sum_{k=1}^Kf(\iota,k;\iota,k)\in\Pc_0. $$
In order to prove (i)--(vi), we will describe a model for part of $p\NvN p$
with respect to the trace $\tau$.
Let $n=K\cdot|E|$.
In the free group factor $L(F_{n^2})$ with tracial state $\tau_F$ consider a
$*$--free family 
$$ S=\bigl((s(\iota,k))_{\iota\in E,\,1\le k\le K},
 (z(\iota_1,k_1;\iota_2,k_2))_{\iota_1,\iota_2\in E,\,1\le k_1,k_2\le K,\,
 (\iota_1,k_1)>(\iota_2,k_2)}\bigr), $$
where with respect to $\tau_F$ each $s(\iota,k)$ is a
semicircular element with second moment $1$ and each
$z(\iota_1,k_1;\iota_2,k_2)$ is a circular element with
$\tau_F\bigl(z(\iota_1,k_1;\iota_2,k_2)^*z(\iota_1,k_1;\iota_2,k_2)\bigr)=1$.
Consider the free product factor
$$ (\RvN,\tau_\RvN)=(L(F_{n^2}),\tau_F)
*\Bigl(\freeprodi(A(\iota),\tau_{A(\iota)})\Bigr), $$
with the embedded copies $A(\iota)\hookrightarrow\RvN$.
Take a system of matrix units
$$ \bigl(g(\iota_1,k_1;\iota_2,k_2)
\bigr)_{\iota_1,\iota_2\in I,\,1\le k_1,k_2\le K} $$
for $M_n(\Cpx)$ and consider the von Neumann algebra $\RvN\otimes M_n(\Cpx)$
with trace $\tau_\RvN\otimes\Tr$, where $\Tr$ is the trace on $M_n(\Cpx)$
taking value $1$ on minimal projections.
Let
$$ \align
Y&=\sum\Sb\iota\in E\\1\le k\le K\endSb s(\iota,k)\otimes g(\iota,k;\iota,k) \\
 \vspace{1.5ex}
&\quad+\sum\Sb\iota_1,\iota_2\in E\\1\le k_1,k_2\le K\\
 (\iota_1,k_1)>(\iota_2,k_2)\endSb
 \Bigl(z(\iota_1,k_1;\iota_2,k_2)\otimes g(\iota_1,k_1;\iota_2,k_2)
 +z(\iota_1,k_1;\iota_2,k_2)^*\otimes g(\iota_2,k_2;\iota_1,k_1)\Bigr) \\
&\hskip30em\in\RvN\otimes M_n(\Cpx).
\endalign $$
Furthermore, consider the subalgebra, $\Bc$, of $\RvN\otimes M_n(\Cpx)$
generated by
$$ \{1_\RvN\otimes g(\iota,j;\iota,k)\mid\iota\in I,\,1\le j,k\le K\}
\cup\bigcup_{\iota\in E}A(\iota)\otimes g(\iota,1;\iota,1). $$
Note that $\Bc$ is isomorphic to $p\MvN p$ in a way that takes the trace
$\tau_\RvN\otimes\Tr$ to $\tau_\MvN$, takes
$A(\iota)\otimes g(\iota,1;\iota,1)$ to $D(\iota)$ and takes
$1_\RvN\otimes g(\iota,j;\iota,k)$ to $e^{(\iota)}_{j,k}$.
Furthermore, by Theorem~2.1 of~\cite{\DykemaRordam}, we have that, with respect
to $n^{-1}(\tau_\RvN\otimes\Tr)$, $Y$ is a semicircular element that is free
from $\Bc$.
{}From the construction of $\NvN$, we have that, with
respect to $\tau(p)^{-1}\tau$, $pXp$ is a semicircular element that is free
from $p\MvN p$.
Hence with $Y$ and $\Bc$ we have constructed a model for $pXp$ and $p\MvN p$.

We identify $\RvN$ with the corner $\RvN\otimes g(\iota_0,1;\iota_0,1)$ of
$\RvN\otimes M_n(\Cpx)$.
For $\iota\in I\backslash\{\iota_0\}$, let
$z(\iota,1;\iota_0,1)=u(\iota)b(\iota)$ be the polar decomposition and let
$u(\iota_0)=1_\RvN$.
Below is a schema indicating the elements (and subalgebras) of $\RvN$
corresponding in the model to the members of the $ F_j^{E,K}$:
$$ \alignat2
 F_1^{E,K}&\ni f(\iota_0,1;\iota,j)Xf(\iota,j;\iota_0,1)
 &&\quad\mapsto\quad u(\iota)^*s(\iota,j)u(\iota) \tag{\Feuone} \\
 F_2^{E,K}&\ni b_\iota=f(\iota_0,1;\iota,1)Xf(\iota_0,1;\iota_0,1)
 &&\quad\mapsto\quad b(\iota)  \tag{\Feutwo} \\
 F_3^{E,K}&\ni f(\iota_0,1;\iota,j)Xf(\iota_0,1;\iota_0,1)
 &&\quad\mapsto\quad u(\iota)^*z(\iota,j;\iota_0,1) \tag{\Feuthr} \\
 F_4^{E,K}&\ni f(\iota_0,1;\iota_1,j)Xf(\iota_2,k;\iota_0,1)
 &&\quad\mapsto\quad u(\iota_1)^*z(\iota_1,j;\iota_2,k)u(\iota_2) \tag{\Feufou} \\
 F_5^{E,k}&\ni u_\iota^*D(\iota)u_\iota
 &&\quad\to\quad u(\iota)^*A(\iota)u(\iota). \tag{\Feufiv}
\endalignat $$
{}From the fact~\cite{\VoiculescuZZCircSemiCirc} that the polar decomposition of a
circular element is of the form $ub$ where $u$ is a Haar unitary, $b$ is a
quarter circular element and $u$ and $b$ are $*$--free, it is easily seen that
each $u(\iota)^*s(\iota,j)u(\iota)$ appearing in~(\Feuone) is a semicircular
element,
each $b(\iota)$ appearing in~(\Feutwo) is a quartercircular
element,
each $u(\iota)^*z(\iota,j;\iota_0,1)$ appearing in~(\Feuthr) is a circular
element and
each $u(\iota_1)^*z(\iota_1,j;\iota_2,k)u(\iota_2)$ appearing in~(\Feufou) is a
circular element.
We must show that the family
$$ \align
 F=\Bigl(\bigl(u(\iota)^*s(\iota,j)u(\iota)&\bigr)_{\iota\in E,
 \,j\in\{1,\ldots,K\}}\,,\,
 \bigl(b(\iota)\bigr)_{\iota\in E\backslash\{\iota_0\}}\,,\,
 \bigl(u(\iota)^*z(\iota,j;\iota_0,1)\bigr)_{\iota\in E,\,j>1}\,,\, \\
&\bigl(u(\iota_1)^*z(\iota_1,j;\iota_2,k)u(\iota_2)\bigr)_{\iota_1,
 \iota_2\in E,\, j,k\in\{1,\ldots,K\},\,
 (\iota_1,j)>(\iota_2,k)>(\iota_0,1)}\Bigr)
\endalign $$
is $*$--free.
Let $z(\iota_1,j;\iota_2,k)=u(\iota_1,j;\iota_2,k)b(\iota_1,j;\iota_2,k)$ be
the polar decomposition.
Taking the polar decomposition of every member of the family $ F$, we see
that this is equivalent to showing that the family
$$ \align
\Bigl(\quad&\bigl(u(\iota)^*s(\iota,j)u(\iota)\bigr)_{\iota\in E,
 \,j\in\{1,\ldots,K\}},\,
 \bigl(b(\iota)\bigr)_{\iota\in E\backslash\{\iota_0\}},\, \\ 
&\bigl(u(\iota)^*u(\iota,j;\iota_0,1)\bigr)_{\iota\in E,\,j>1},\, 
 \bigl(b(\iota,j;\iota_0,1)\bigr)_{\iota\in E,\,j>1},\, \\
&\bigl(u(\iota_1)^*u(\iota_1,j;\iota_2,k)u(\iota_2)\bigr)_{\iota_1,
 \iota_2\in E,\, j,k\in\{1,\ldots,K\},\,
 (\iota_1,j)>(\iota_2,k)>(\iota_0,1)}, \\
&\bigl(u(\iota_2)^*b(\iota_1,j;\iota_2,k)u(\iota_2)\bigr)_{\iota_1,
 \iota_2\in E,\, j,k\in\{1,\ldots,K\},\,
 (\iota_1,j)>(\iota_2,k)>(\iota_0,1)}\quad\Bigr)
\endalign $$
is $*$--free.
However, replacing every quartercircular $b(\cdots)$ and every semicircular
$s(\cdot,\cdot)$ by Haar unitaries that generate the same algebra, the
$*$--freeness of~(\bigffam) reduces to a question of freeness of a certain set
of words in a free group; these words are  seen by inspection to be
free, because each of the words in question has one letter that none of the
others has.

Let $\Ft$ be the set of all the elements that appear in the family $ F$.
To show that $\QvN\cong A$, it remains only to show that the family
$$ \Bigl(\Ft,\bigl(u(\iota)^*A(\iota)u(\iota)\bigr)_{\iota\in E}\Bigr)
\tag{\Feutfam} $$
of sets of random variables is $*$--free.
This is easily shown by induction on $|E|$.
Indeed, let $Z$ be the $*$--algebra generated by the polar and positive parts
of all the $z(\cdot,\cdot;\cdot,\cdot)$'s.
Then $\Ft\subseteq Z$, each $u(\iota)\in Z$ and, by assumption, the family
$\bigl(Z,(A(\iota))_{\iota\in E}\bigr)$ is free.
Let $\iota_1\in E$;
then
$\bigl(Z,(u(\iota)^*A(\iota)u(\iota))_{\iota\in E\backslash\{\iota_1\}}\bigr)$
is free by inductive hypothesis (or tautologically if $|E|=1$).
Hence to show that~(\Feutfam) is free, it will suffice to show that
$\Zt\eqdef\alg\bigl(Z\cup\bigcup_{\iota\in E\backslash\{\iota_1\}}A(\iota)\bigr)$ and
$u(\iota_1)^*A(\iota_1)u(\iota_1)$ are free.
But since $u(\iota_1)$ belongs to $\Zt$, after conjugating by $u(\iota_1)$ we
simply appeal to freeness of $A(\iota_1)$ and $\Zt$.
This completes the proof that $\QvN\cong A$.

Suppose now that each of the II$_\infty$ factors $A(\iota)\otimes B(\Hil)$ has
a continuous one--parameter group of trace--scaling automorphisms.
Taking the direct sum of these automorphisms, we get a continuous
one--parameter group, $(\alpha_\lambda)_{\lambda\in\Reals_+^*}$ of
automorphisms of $\MvN$ such that
$\tau_\MvN\circ\alpha_\lambda=\lambda\tau_\MvN$.
We will construct a continuous one--parameter group of automorphisms
$(\gamma_\lambda)_{\lambda\in\Reals_+^*}$ of $\NvN$ such that $\tau\circ\gamma_\lambda=\lambda\tau$.
Let $\beta_\lambda:\Cpx[X]\to\Cpx[X]$ be the automorphism defined by
$\beta_\lambda(X)=\lambda^{-1/2}X$;
given $p\in\Pc_0$ let
$\gamma_\lambda=\alpha_\lambda*\beta_\lambda:\Ac_p\to\Ac_{\alpha_\lambda(p)}$
be the resulting isomorphism.
Using again Speicher's combinatorial approach to freeness and the
expression~(\psipexpr), and recalling that $\pi^c$ has $1+n/2$ blocks, we get
$$ \align
\psi\circ&\gamma_\lambda(a_0Xa_1Xa_2\cdots Xa_n)= \\
&=\lambda^{-n/2}\psi(\alpha_\lambda
 (a_0)X\alpha_\lambda(a_1)X\alpha_\lambda(a_2)\cdots X\alpha_\lambda(a_n))
 \\
&=\lambda^{-n/2}\sum_{\pi\in\NCP(\{2,4,\ldots,2n\})}(\phi)_{\pi^c}
 [\alpha_\lambda(a_0),\alpha_\lambda(a_1),\alpha_\lambda(a_2),\ldots,\alpha_\lambda(a_n)] \\
&=\lambda\sum_{\pi\in\NCP(\{2,4,\ldots,2n\})}(\phi)_{\pi^c}
 [a_0,a_1,a_2,\ldots,a_n] \\
&=\lambda\psi(a_0Xa_1Xa_2\cdots Xa_n).
\endalign $$
Thus $\psi\circ\gamma_\lambda=\lambda\psi$ on $\Ac$;
using $\lambda^{-1/2}\gamma_\lambda$ to define a unitary operator $U_\lambda$
on $L^2(\Ac,\psi)$, we see that $\gamma_\lambda$ extends via conjugation by
$U_\lambda^*$ to an automorphism of $\NvN$, which we also call
$\gamma_\lambda$.
Since the unitary group $(U_\lambda)_{\lambda\in\Reals_+^*}$ is strongly
continuous, the group $(\gamma_\lambda)_{\lambda\in\Reals_+^*}$ is $\sigma^*$--strongly
continuous.
\QED

\proclaim{Corollary
 \IIIone}
For every $n\in\Nats$ let $A(n)$ be a II$_1$ factor such that the associated
II$_\infty$ factor $A(n)\otimes B(\Hil)$ has a continuous one--parameter
group of trace--scaling automorphisms.
Let $\MvN=\freeprod{n=1}\infty A(n)$ be their free product (with respect to
tracial states).
Then $\MvN\otimes B(\Hil)$ has a continuous one--parameter group of
trace--scaling automorphisms, and hence there is a type III$_1$ factor whose
core is $\MvN\otimes B(\Hil)$.
\endproclaim
\demo{Proof}
An application of Theorem~\IIinf{} gives that
$\bigl(\MvN*L(F_\infty)\bigr)\otimes B(\Hil)$ has a continuous one parameter
group of trace--scaling automorphisms.
But Theorem~\Compr{} shows that $\MvN*L(F_\infty)\cong\MvN$.
\QED

\enddocument